\documentclass[journal, 11pt]{IEEEtran}

\ifCLASSINFOpdf

\else

\fi
 
\usepackage{url}
\usepackage{mathtools, cuted}
\usepackage{stfloats}
\usepackage{hyperref}
\usepackage{amsmath,amsfonts,amssymb}
\usepackage{wrapfig}
\usepackage{draftwatermark}
\SetWatermarkText{}
\SetWatermarkScale{1}
\newtheorem{definition}{Definition}

\newtheorem{lemma}{Lemma}
\newtheorem{theorem}{Theorem}
\newtheorem{comment}{Comment}
\newtheorem{corollary}{Corollary}
\newtheorem{example}{Example}
\usepackage{mathptmx} 
\usepackage{amsmath}  
\usepackage{amssymb}  
\usepackage{amsfonts}
\usepackage{boxedminipage}
\usepackage{graphicx}
\usepackage[T1]{fontenc}
\usepackage[english]{babel}
\usepackage{pst-sigsys}
\usepackage{pst-node,pst-circ}
\definecolor{dgray}{gray}{0.6}
\definecolor{lgray}{gray}{0.8}
\newcommand{\ignore}[1]{}

\begin{document}

\title{Stability tests for a class of switched descriptor systems with non-homogenous indices}


\author{Shravan~Sajja,~Martin~Corless, Ezra Zeheb,~Robert Shorten
\thanks{S. Sajja  is with IBM Research Ireland.}
\thanks{M. Corless is with the School of Aeronautics and Astronautics, Purdue University, West Lafayette, IN, USA.}
\thanks{E. Zeheb is with the Technion-Israel Institute of Technology, Haifa and Holon College of Engineering, Holon, Israel.}
\thanks{R. Shorten is with University College Dublin, Ireland. }}

\maketitle

\parindent0mm

\begin{abstract}

In this paper we derive stability conditions for a switched system where  switching occurs between linear  descriptor systems of different indices.
In particular, our results can be used to analyse the stability of the important case when switching between a standard system and an index one descriptor system, and 
systems where switching occurs between  an index one and an index two descriptor system.  Examples are given to illustrate the use of our results. 
\end{abstract}


\IEEEpeerreviewmaketitle

\section{Introduction}

\noindent 

Descriptor systems provide a natural framework to model and analyse many dynamic systems with algebraic constraints. 
They appear frequently in modelling engineering systems: for example in the description of interconnected
large scale systems; in economic systems (e.g.
the fundamental dynamic Leontief model); network analysis \cite{dai1989singular} and they
are also particularly important in the simulation and design
of very large scale integrated (VLSI) circuits.\newline

Recently, motivated by certain applications, some authors have begun the study of descriptor systems that are characterized
by switching between a number of descriptor modes \cite{s2009trenn,liberzon2009stability,zhai2011commutation, zhou2013, wirth2015, {strenn2012}, {trenn2012}}. 
For example, in \cite{liberzon2009stability} the authors focus on dwell time arguments, and on conditions on the "consistency projectors", to obtain stability under arbitrary
switching. In \cite{wirth2011commutation} and  \cite{zhai2011commutation}, the authors, under an assumption of a state-dependent switching
condition (to avoid impulses), obtain a condition for stability based on commuting vector fields. These results mimic similar results derived for standard switched systems by \cite{mori1997, morse1999liberzon, shorten2001}.  \newline

Our approach in this paper differs from that given in the above papers. Our results in this note are based on a fundamental result derived in \cite{dimensionality} to recursively reduce the dimension of a switched descriptor system using full rank decomposition.  This approach allows us to obtain conditions
which can also be checked without resorting to complicated linear algebraic manipulations. 
This has been achieved for a special system class of   index one descriptor systems; namely, a class of switched systems characterised by rank-1 perturbations,
for which a simple continuity assumption on the state at the switching instances is satisfied (see \cite{dimensionality}). In this note, we now extend the results in \cite{dimensionality} to switching between an index one descriptor system and (i)
a standard system (which can be described completely by a set of ordinary differential equations); and (ii) an index  two descriptor system.  some results presented earlier in \cite{DescripConf2013}. The present results provide detailed proofs and explanations for some results presented earlier in \cite{DescripConf2013}, and provide new examples of systems to which these results can be applied. In particular, our results apply to systems for which standard assumptions in the descriptor literature, do not apply.

\section{Preliminary Results}

\noindent Consider a  linear time invariant (LTI) descriptor system described by
\begin{equation}\label{adescrip}
E\dot{x}=Ax,
\end{equation}
where $ E, A \in \mathbb{R}^{n\times n}$. When $E$ is nonsingular, this system is also described by the {\it standard system} $\dot{x} = E^{-1}Ax$.
When $E$ is singular, then both algebraic equations and differential equations describe the behavior of the system, and the system is known 
as a {\em descriptor system.}
Since we shall be interested in switched systems which are constructed by switching 
between systems that are  exponentially stable about zero, we require  $A$ to be nonsingular \cite{owens1985};
note that were $A$ singular,  there would be equilibrium states other than zero.\newline

 The following notions are important when studying descriptor systems.
  First, the system \eqref{adescrip} is said to be \underline{\it stable} if  every   eigenvalue of 
  $(E,A)$ has a negative real part and the system is called \underline{\it regular} if $\det[sE-A]\not\equiv 0$.\newline 

%
%
With $A$ invertible, the \underline{\em index} of system \eqref{adescrip} or the pair $(E, A)$  is the smallest integer $k^* \leq n$ for which
\begin{equation}
\label{eq:index}
\text{Im}((A^{-1}E)^{k^*+1}) = \text{Im}((A^{-1}E)^{k^*})
\end{equation}
where Im denotes the image of a matrix. Thus,  a standard  system is  an index zero descriptor system.\newline  

The \underline{\it consistency space} for system (\ref{adescrip}) or $(E,A)$ is defined by
\begin{equation}
\mathcal{C}= \mathcal{C}(E,A) := \text{Im}((A^{-1}E)^{k^{*}})
\end{equation}
where $k^*$ is the index of $(E,A)$.
 Note that $\mathcal{C}$ is the set of all initial states  for which the system has a continuous solution.
Since $\text{Im}((A^{-1}E)^{k^*+1}) = \text{Im}((A^{-1}E)^{k^*})$ we see that $A^{-1}E \mathcal{C} = \mathcal{C}$;
this means that $A^{-1}E$ is a one-to-one mapping of $\mathcal{C}$ onto itself; hence the kernel of $E$
and $\mathcal{C}$ intersect only at zero \cite{owens1985}. 
If $\mathcal{C}=\{0\}$ the system is trivial and the only continuous solution is the zero solution $x(t) \equiv 0$.
If $\mathcal{C} \neq \{0\}$,  we let $\tilde{A}$ be the inverse of the map $A^{-1}E$ restricted to
$\mathcal{C}$; then (\ref{adescrip}), or equivalently $x=A^{-1}E\dot{x}$, is equivalent to
\begin{equation}
\label{eq:descriptor3}
\dot{x} = \tilde{A}x
\end{equation}
Thus the  restriction of the descriptor system to its consistency space is equivalent to the standard system (\ref{eq:descriptor3}) where $x(t)$ is in $\mathcal{C}$.\newline

\noindent Another way to introduce the consistency space is as follows. When $(E, A)$ is regular with a non-trivial consistency space,
 it can be shown
   that there exist nonsingular matrices $S$ and $T$ such
that  \cite{trenn2009}
\[
SET = \left[
\begin{array}{cc} I	& 0	\\
0	&N
\end{array}\right]
\qquad 
SAT = \left[
\begin{array}{cc} J	& 0	\\
0	&I
\end{array}\right]
\]
where the matrix $N$ is nilpotent, that is,  $N^k =0$ for some $k\ge 1$.
If $A$ is nonsingular, then $J$ is nonsingular and
for any $k \ge 1$,
\[
(A^{-1}E)^k = 
T \left[
\begin{array}{cc} J^{-k}	& 0	\\
0	&N^k
\end{array}\right]
T^{-1}
\]
Then the index of $(E, A)$ is the smallest $k^*$ for which $N^{k^*} =0$ and
the consistency space is the range of the consistency projector defined by
\[
\Pi = T \left[
\begin{array}{cc} I	& 0	\\
0	&0
\end{array}\right]
T^{-1}
\]
Also
\[
\tilde{A}= T \left[
\begin{array}{cc} J	& 0	\\
0	&0
\end{array}\right]
T^{-1}
\]
Note that the  consistency space  for $(E, A)$ is trivial if and only if $A^{-1}E$ is nilpotent.

\subsubsection*{Lyapunov functions}
To obtain stability conditions for switching descriptor systems we shall use Lyapunov functions.
Let us first consider Lyapunov functions for system \eqref{adescrip}.
Consider any differentiable function $V:\mathbb{R}^n \rightarrow \mathbb{R}$.
Its derivative  along solutions of  \eqref{adescrip} is given by $\dot{V} = DV(x)\dot{x}$ which can be expressed as a function of $x$. 

\begin{definition}
A  differentiable function $V:\mathbb{R}^n \rightarrow \mathbb{R}$ is a \underline{\it Lyapunov function} for  $(E, A)$ if $V$ is positive definite on 
$\mathcal{C}=\mathcal{C}(E, A)$
and $\dot{V}$ is negative for any non-zero state in $\mathcal{C}(E, A)$.
A symmetric matrix $P$ is a \underline{\it Lyapunov matrix} for $(E, A)$  if $V(x) = x^TPx$ is a Lyapunov function for $(E, A)$.\newline
\end{definition}

The following lemma provides a characterization of all Lyapunov matrices for a linear descriptor system.\newline

\begin{lemma}
\label{lemma:Lyap}
A symmetric matrix  $P$  is a Lyapunov matrix for $(E, A)$ if and only if
$P$   is positive-definite  on
$\mathcal{C}=\mathcal{C}(E, A)$
and 
$PA^{-1}E+ E^TA^{-T}P$
is negative definite
 on $\mathcal{C}$.
\end{lemma}

 \noindent{\sc Proof:} 
It follows  from  \eqref{adescrip}  and the invertibility of $A$ that $x=A^{-1}E\dot{x}$; hence
\begin{align*}
\dot{V}= 2x^TP\dot{x} &
=2(A^{-1}E\dot{x})^TP\dot{x}  = \dot{x}^T(PA^{-1}E+E^TA^{-T}P)\dot{x}\\
&= -\dot{x}^TQ\dot{x} 
\end{align*}
where
$
Q:=
-PA^{-1}E -E^TA^{-T}P
$.
Recall that system description (\ref{adescrip}) is equivalent to $\dot{x}= \tilde{A}x$,
where $x$ is in $\mathcal{C}$
and $\tilde{A}$ is an invertible map on $\mathcal{C}$. 
Hence
\[
\dot{V} =  -x^T\tilde{Q}x \quad \mbox{where} \quad \tilde{Q} = -\tilde{A}^TQ\tilde{A}.
\]
Since $\tilde{A}$ maps $\mathcal{C}$ onto $\mathcal{C}$ and is invertible on $\mathcal{C}$,
  $\tilde{Q}$  is positive-definite on $\mathcal{C}$ if and only if
$Q$ is positive-definite on  $\mathcal{C}$.
Hence, $\dot{V}$   is negative for any non-zero state in $\mathcal{C}$ if and only if $PA^{-1}E +E^TA^{-T}P$ is negative-definite on $\mathcal{C}$.\textbf{ Q.E.D.}
\newline


Previous papers  such as  \cite{owens1985,trenn2009,liberzon2009stability,stykelPhD2002} consider  a specific class of 
Lyapunov matrices of the form $P=E^T\tilde{P}E$ where $\tilde{P}$ is a positive definite matrix for which
$E^T\tilde{P}A +A^T\tilde{P}E$ is negative definite on $\mathcal{C}$. 
In particular, \cite{owens1985} shows that the existence of a Lyapunov matrix of this type is necessary and sufficient for asymptotic stability of system \eqref{adescrip}.\newline

Let $C$ be any matrix with the following property:
\begin{equation}
\label{eq:C}
x\in \mathcal{C}(E,A) \quad \mbox{if and only if} \quad Cx = 0
\end{equation}
Then we have the following LMI characterization of Lyapunov matrices for $(E,A)$.

\begin{lemma}
\label{lem:LMI}
A symmetric matrix $P$ is a Lyapunov matrix for $(E, A)$ if and only if there exists  scalars $\kappa_1, \kappa_2 \ge 0 $ such that
\begin{align}
 P + \kappa_1C^TC& >0		\label{eq:LMI0}\\
 PA^{-1}E +E^TA^{-T}P -\kappa_2 C^TC  &<0			\label{eq:LMI1}
\end{align}
where $C$ is any    matrix satisfying \eqref{eq:C}.
\end{lemma}

{\sc Proof:}
A vector $x$ is in $\mathcal{C} = \mathcal{C}(E,A)$ if and only if   $Cx = 0$ which is equivalent to
$x^TC^TCx = 0$. From Lemma \ref{lemma:Lyap} we know  that $P$ is a Lyapunov matrix for $(E, A)$ if and only if $x^TPx >0$ and
$x^T(PA^{-1}E +E^TA^{-T}P)x< 0$ whenever $x\neq 0$ and $x^TC^TCx = 0$.
It follows from Finslers Lemma that there exist scalars $\kappa_1$ and $\kappa_2$ such that \eqref{eq:LMI0} and  \eqref{eq:LMI1} hold;
since $C^TC \ge 0$, \eqref{eq:LMI0}  and  \eqref{eq:LMI1} also hold with $\kappa_1, \kappa_2 \ge 0$.

Note that if there exists a Lyapunov matrix $P$ satisfying \eqref{eq:LMI0}  and  \eqref{eq:LMI1} then there also exists a Lyapunov matrix satisfying

\begin{align}
\label{eq:LMI2}
 P + C^TC &>0\\
 PA^{-1}E +E^TA^{-T}P - C^TC  &<0
\end{align}


\subsection{Linear switched descriptor systems } 
\noindent The ultimate objective of our work is to analyze the stability
of  switched descriptor systems described by
\begin{equation}
\label{eq:descriptor0}
E_{\sigma(t)} \dot x = A_{\sigma(t)}x
\qquad \mbox{where} \qquad
 \sigma(t) \in \{  1, \cdots , N\}
\,.
\end{equation}
We assume throughout this paper that $\sigma$ is  piecewise continuous  with a finite number of discontinuities
in any bounded time interval. Also,  {\em we  do not  require \eqref{eq:descriptor0} to hold at points of discontinuity of $\sigma$.}
Thus, if $\sigma$ is continuous at $t$ and $\sigma(t)=i$,
the system must satisfy
\[
E_i \dot x = A_i x\,;
\]
hence $x(t)$ must be in the consistency space of $(E_i, A_i)$.
To complete the description of the  switching descriptor systems under consideration one must  also specify
how the system behaves
during switching.
If $\sigma$  switches from $i$ to $j$ at $t_*$ then $x(t_*^-):=\lim_{t \rightarrow t_*, t < t_*}x(t)$ must be
in $\mathcal{C}(E_i, A_i)$ and $x(t_*^+):=\lim_{t \rightarrow t_*, t > t_*}x(t)$ must be
in $\mathcal{C}(E_j, A_j)$.
If $x(t_*^-)$ is not in  $\mathcal{C}(E_j, A_j)$ then, one has to have a solution
which is discontinuous at $t_*$ and  to complete the description one must specify how $x(t_*^+)$ is obtained from $x(t_*^-)$.
In some switching systems, there is a restriction on $x(t_*^-)$, that is, $\sigma$ can only switch from a specific $j$ to a specific $i$ if 
$x(t_*^-)$ is in a  restricted  subset of the consistency space of $(E_i, A_i)$.
This is illustrated in  Example \ref{mechanical}.\newline

In some treatments of switched descriptor systems, 
  \eqref{eq:descriptor0}  is satisfied for all $t$
\cite{trenn2009,liberzon2009stability}.
In that case, one has to consider $x(\cdot)$ to be a distribution because the solution  of  \eqref{eq:descriptor0}  may contain impulses.
When \eqref{eq:descriptor0} is  satisfied for all $t$ and $\sigma$ is continuous from the right, it can be shown   that,  if $\sigma$ switches from
$i$ to $j$ at $t_*$ then
\begin{equation}
\label{eq:jump}
x(t_*^+)= \Pi_j x(t_*^-)
\end{equation}
where
$\Pi_j$ is the consistency projector associated with $(E_j, A_j)$.
Here we do not require that \eqref{eq:descriptor0} be satisfied at discontinuity points of $\sigma$ nor do we require 
\eqref{eq:jump}.
This is illustrated in  Example \ref{mechanical}.
When a  system satisfies \eqref{eq:descriptor0} for all $t$
 and if switching can occur from any state in
the consistency space of $(E_i, A_i)$ then, as in \cite{liberzon2009stability} one must assume that
\begin{equation}
E_j(I-\Pi_j)\Pi_i = 0
\end{equation}
in order to guarantee solutions without impulses  when switching from  from
$i$ to $j$. We do not need this assumption here.

%



\subsection{Lyapunov stability conditions}

\noindent Our next result contains
conditions
that are sufficient to guarantee  \underline{\it global uniform exponential stability (GUES)}  of \eqref{eq:descriptor0}.\newline

\begin{theorem}
\label{lemma:qs}
Consider a switched descriptor system which satisfies
(\ref{eq:descriptor0})
at points of continuities of $\sigma$ and  suppose that  for each  $i=1, \cdots, N$,
  there is a Lyapunov  matrix $P_i$  for $(E_i, A_i)$ 
such that
 \begin{equation}\label{lyap_inequality}
 x(t_*^+)^TP_jx(t_*^+) \le x(t_*^-)^TP_ix(t_*^-)\,.
 \end{equation}
whenever    $\sigma$ switches from $i$ to $j$ at $t_*$.
Then the system is {GUES}.
\end{theorem}
\vspace{0.5em}

 \noindent{\sc Proof:} 
Consider any solution $x(\cdot)$ of the system, and let $v(t) =x(t)^TP_{\sigma(t)}x(t).$
If $t$ is a point of discontinuity of $\sigma$, then, according to \eqref{lyap_inequality},
\begin{equation}
\label{eq:Vdisc}
v(t^+) \le v(t^-).
\end{equation}
Suppose  $t$ is not a point of discontinuity of $\sigma$.
If $\mathcal{C}_i =\mathcal{C}(E_i, A_i)= \{0\}$ then $v(t) = 0$ and  $\dot{v}(t) = 0$.
Otherwise, 
$
 E_i\dot{x}=A_ix
$
  where $i =\sigma(t)$.  
 Following the proof of Lemma \ref{lemma:Lyap},  
\[
\dot{v} = -x^T\tilde{Q}_ix \quad \mbox{where} \quad \tilde{Q}_i = \tilde{A}_i^TQ_i\tilde{A}_i
\]
 \begin{equation}\label{def_Qi2}
Q_i  = -P_iA_i^{-1}E_i -E_i^TA_i^{-T}P_i
\end{equation}
and $\tilde{Q}_i$ is positive-definite on  $\mathcal{C}_i$.
Recalling that $P_i$ is positive-definite on $\mathcal{C}_i$, let
\[
\alpha_i = \frac{1}{2}\min\{x^T\tilde{Q}_i x: x \in \mathcal{C}_i  \mbox{ and } x^TP_ix = 1 \}.
\]
Then $\alpha_i >0$ and $\dot{v} \le - 2 \alpha_i v$.
Now let $\alpha = \min\{ \alpha_1, \cdots \alpha_n \}$.
Then $\alpha >0$ and
\begin{equation}
\dot{v}(t) \le -2 \alpha v(t)
\end{equation}
when $\sigma$ is continuous at $t$. From this and the discontinuity condition (\ref{eq:Vdisc}), we can conclude that
$ v(t) \le e^{-2\alpha (t-t_0)}v(t_0)$ for $t \ge t_0$.
Since each $P_i$ is positive-definite on $\mathcal{C}_i$
  there are constants $\lambda_1, \lambda_2 >0$ such that, for $i=1, \cdots, N$, we have 
$\lambda_1\|x\|^2 \le  x^TP_ix \le \lambda_2\|x\|^2$ whenever $x$ is in  $\mathcal{C}_i$.
 Hence
$
\lambda_1 \|x(t)\|^2\le  \nu(t) \le \lambda_2 \|x(t)\|^2
$
 and every solution $x(\cdot)$ satisfies
\begin{equation}
\|x(t)\| \le \beta e^{-\alpha (t-t_0)}\|x(t_0)\|
\end{equation}
 for all $ t\ge t_0$,
  where $\beta=\sqrt{\lambda_2/\lambda_1}$.
This means that the system is GUES.
\textbf{Q.E.D.}\newline
 
We have the following corollary to Theorem \ref{lemma:qs} and Lemma \ref{lem:LMI}.

\begin{corollary}\label{lemma:cor1}
Consider a switching descriptor system described by
(\ref{eq:descriptor0})
and suppose that
 there is a symmetric   matrix $P$ satisfying
\begin{align}
P +C_i^TC_i >0\\
PA_i^{-1}E_i + E_i^TA_i^{-1}P  -C_i^TC_i  <0
\end{align}
for $i=i, \dots, N$ where
\begin{equation}
\label{eq:Ci}
x\in \mathcal{C}(E_i,A_i) \quad \mbox{if and only if} \quad C_ix = 0
\end{equation}
Also,
\begin{equation}
\label{eq:vdecrease}
x(t_*^+)Px(t_*^+) \le x(t_*^-)^TPx(t_*^-)
\end{equation}
if $\sigma$ switches at $t_*$.
 Then, the system is {GUES}.
 \newline
 \end{corollary}

To conclude this section, we present an example to motivate our results. 
This example illustrates the use of Theorem \ref{lemma:qs} to analyse stability of 
switching between a standard  system and a descriptor system.\newline

\begin{example}[A simple switched mechanical system]
\label{mechanical}
\begin{figure}[h]
\begin{center} \includegraphics[height=1in]{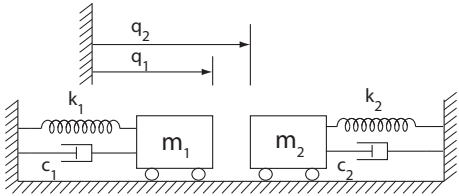}
\caption{A mechanical system}
\end{center}
\label{fig:mech_sys1}
\end{figure}
\noindent Consider the  mechanical system  illustrated in Figure \ref{fig:mech_sys1}
consisting of two spring mass dampers with
masses $m_1, m_2$, damping coefficients $c_1, c_2$ and spring constants $k_1, k_2$, respectively.
Let $q_1$ and $q_2$ denote the displacements of the masses $m_1$ and $m_2$ from their rest
positions.
We consider the switched system in which the two masses can lock onto each other
when their displacements are the same; see Figure \ref{fig:mech_sys2}.
When they are locked together their displacements remain equal.
When they unlock, their displacements are independent.

\begin{figure}[h]
\centering
\includegraphics[height=1in]{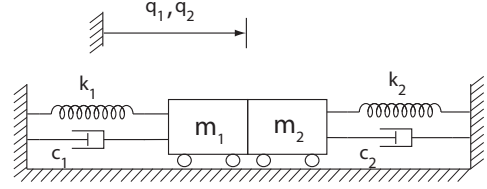}
\caption{Masses locked together}
\label{fig:mech_sys2}
\end{figure}

When the masses are not locked together, the system is described by
\begin{equation}
\begin{array}{rcl}
m_1\ddot{q}_1 + c_1 \dot{q}_1 +k_1q_1 &=&0\\
m_2\ddot{q}_2 + c_2 \dot{q}_2 +k_2q_2 &=& 0
\end{array}
\end{equation}
When the masses are locked together, we have the following description
\begin{equation}
\begin{array}{rcl}
m\ddot{q}_1 + c\dot{q}_1 +kq_1 &=& 0\\
q_1 -q_2 &=& 0
\end{array}
\end{equation}
where
\[
m:=m_1+m_2\,, \quad c:=c_1+c_2\,, \quad k:=k_1 +k_2 \,.
\]
Since lock-up is due to internal forces in the system,
linear momentum is conserved during lock-up, that is,
if lockup occurs at time $t$ then,
\[
m\dot{q}_2(t^+) =m\dot{q}_1(t^+) = m_1\dot{q}_1(t^-) + m_2 \dot{q}_2(t^-)
\]
which results in
\begin{equation}
\dot{q}_2(t^+) =\dot{q}_1(t^+)= \frac{m_1\dot{q}_1(t^-) +m_2 \dot{q}_2(t^-)}{m}
\end{equation}
We also have
\begin{eqnarray}
q_1(t^+) =q_2(t^+)= q_1(t^-)= q_2(t^-)
\end{eqnarray}

\noindent Introducing state variables $x_1 = q_1, x_2 = \dot{q}_1, x_3 = q_2, x_4 = \dot{q}_2$,
this system can be described by the switched system (\ref{eq:descriptor0}) where $N=2$,
{\small
\begin{align*}
&E_1 =
\left[
\begin{array}{cccc}
1   &0  &0  &0\\
0   &m_1  &0  &0\\
0   &0  &1 &0\\
0   &0  &0  &m_2
\end{array}
\right]
\,\,\,
 A_1 =
\left[
\begin{array}{rrrr}
0   &1  &0  &0\\
-k_1   &-c_1  &0  &0\\
0  &0  &0 &1\\
0     &0    &-k_2  &-c_2
\end{array}
\right]
\\
&E_2 =
\left[
\begin{array}{cccc}
1   &0  &0  &0\\
0   &m &0  &0\\
0   &0  &1 &0\\
0   &0  &0  &0
\end{array}
\right]
\,\,\,
\quad A_2 =
\left[
\begin{array}{rrrr}
0   &1  &0  &0\\
-k   &-c  &0  &0\\
  0 &0   &0 &1\\
-1  &0  &1 &0
\end{array}
\right]
\end{align*}}
Observe that $(E_1, A_1)$ is a standard system  and it can be readily shown that  $(E_2, A_2)$ is an index two descriptor system. 

 \noindent
 Switching    to mode one  can  occur  at any state in the consistency space of mode two and if this occurs at time $t_*$, then
 $x(t_*^+) = x(t_*^-)$.
 Switching   to mode two  only occurs when $x_1=x_3$ and if this occurs at time $t_*$,
\begin{equation}
\label{eq:switch_x}
\begin{array}{rcl}
x_1(t_*^+) = x_3(t_*^+) &=& x_1(t_*^-)=x_3(t_*^-)\\
x_2(t_*^+) = x_4(t_*^+) &=&  \frac{m_1 x_2(t_*^-) +  m_2 x_4(t_*^-)}{m}
\end{array}
\end{equation}
Note that this system does not satisfy condition \eqref{eq:jump} at switching. One can show that this condition requires that
$x_2(t_*^+) = x_4(t_*^+) = x_2(t_*^-)$.

\noindent As  candidate Lyapunov matrices for this system consider
\[
P_1 = P_2 = P
: =
\left[
\begin{array}{cccc}
k_1  &\epsilon m_1  &0  &0\\
\epsilon m_1   &m_1  &0  &0\\
0   &0  &k_2 &\epsilon m_2\\
0   &0  &\epsilon m_2  &m_2
\end{array}
\right]
\]

\noindent where $\epsilon >0$; clearly $P>0$ for $\epsilon$ sufficiently  small.
Recalling definition (\ref{def_Qi2}) of  the symmetric matrix $Q_i$ we obtain
{\small
\[
Q_1 =
\left[
\begin{array}{cccc}
2c_1\!-\!2\epsilon m_1               &\epsilon \frac{m_1c_1}{k_1}           &0  &0\\
*  &2\epsilon \frac{m_1^2}{k_1} &0  &0\\
0  &0  &2c_2\!-\!2\epsilon m_2  &\epsilon \frac{m_2c_2}{k_2}\\
0     &0    &*  &2\epsilon \frac{m_2^2}{k_2}
\end{array}
\right] \]
\[
Q_2 =
\frac{1}{k}\left[
\begin{array}{cccc}
2k_1c\!-\!2\epsilon m_1k 				&  mk_1\!-\!m_1k   \!+\! \epsilon m_1c 	&k_2  c  & \epsilon m_2 c \\
*								&2\epsilon m m_1 				&mk_2  &\epsilon m m_2\\
*			 					&*							&-2\epsilon_2 k  &-m_2k\\
* 								& *  							&* &0
\end{array}
\right]
\]
}Clearly, $Q_1 >0$ for $\epsilon$ sufficiently small.
When $x$ is in the consistency space of $(E_2, A_2)$,
we have $x_3 = x_1$, $x_4 = x_2$; hence
\[
x^TQ_2x =
\left[\begin{array}{c}x_1\\x_2 \end{array}\right]^T
\tilde{Q}_2
\left[\begin{array}{c}x_1\\x_2 \end{array}\right]^T\]
\noindent where
\[\tilde{Q}_2=
\frac{1}{k}\left[
\begin{array}{cc}
2kc -2\epsilon m k &     \epsilon mc  \\
    \epsilon mc & 2\epsilon m^2
\end{array}
\right]
\]
and $\tilde{Q}_2 >0$ for $\epsilon$ sufficiently small; hence $Q_2$ is positive definite on the consistency space of $(E_2, A_2)$.

\noindent
When switching   to mode one,
$x(t_*^+)^TP_1x(t_*^+) =x(t_*^-)^TP_2x(t_*^-) $.
When switching  to mode two, it follows from    the switching conditions in (\ref{eq:switch_x})  that
{\small
\[
x(t_*^+)^TP_2x(t_*^+) - x(t_*^-)^TP_1x(t_*^-) =
-\frac{m_1m_2}{m}[x_2(t_*^-)-x_4(t_*^-)]^2  
\le 0
\]
}
\noindent
Thus, for $\epsilon$ sufficiently small, the matrices $P_1$  and $P_2$ satisfy the requirements of Lemma \ref{lemma:qs};
hence this system is globally uniformly exponentially stable about zero for any allowable switching sequence.
\end{example}


\section{Main results}
\label{sec:indexone}


We now present some  simple tests to check the stability of special classes of switched descriptor systems  constructed by
(i) switching between a standard  system and an index-one system; and (ii) switching between index-one and index-two descriptor systems.
These results build on  \cite{shorten2009quadratic} and the following lemma.

\vspace{1em}
\begin{lemma}
\label{lem:interlacing}
Suppose $A, B  \in \mathbb{R}^{n \times n}$ with  $rank(A-B)= rank(A) -rank(B)$ and
\begin{eqnarray*}
A+ A^T  &>&   0	\\
 B  + B^T&\geq&  0.
\end{eqnarray*}
Then,
the kernels of $B$ and $B+B^T$ are equal.
\end{lemma}
\vspace{1em}

\noindent{\sc Proof:} Since $Q_A >0$, where $Q_A:=A+A^T$, we see that $rank(A) = n$.
Let $r:= rank(B)$. Then, by assumption, we have
$rank(A-B) = n-r$.
 Recall that the nullity of a matrix is the dimension of its kernel.
First, we show that
the nullity of $Q_B:=B+B^T$ is at most $n\!-\!r$.
So, suppose that $x\neq 0$ is in the kernel of $Q_B$. Then
\begin{eqnarray}
0  &=&x^TQ_Bx =x^T(A+A^T)x + 2x^T(B\!-\!A)x \nonumber \\ &=&  x^TQ_Ax -2x^T(A\!-\!B)x \nonumber
\end{eqnarray}
Since $Q_A > 0$ and $x\neq 0$, we have $x^TQ_Ax >0$;
hence   $(A\!-\!B)x \neq 0$, that is, $x$ is not in the kernel of $A\!-\!B$.
Thus, the kernel of $Q_B$ and $A\!-\!B$ intersect only at zero.
Since the rank of $A\!-\!B$ is $n\!-\!r$, its nullity is $r$; hence the nullity of $Q_B$ is at most $n\!-\!r$.
We now show that the kernel of $Q_B$ contains the kernel of $B$.
So, suppose  that $x$ is in the kernel of $B$, that is, $Bx=0$. Then
$
x^TQ_Bx = 2x^TBx = 0
$.
Since $Q_B \ge 0$, it follows that $Q_Bx = 0$, that is, $x$ is in the kernel of $Q_B$.
Thus, the kernel of $Q_B$ contains the kernel of $B$. 
Finally, since $B$ has rank $r$, its nullity is
 $n\!-\!r$. Since we also know that the nullity  of $Q_B$ is less than or equal to $n\!-\!r$,
it now follows that the kernel of $Q_B$ is the same as the kernel of $B$. \textbf{Q.E.D.}\newline

The following general result is  a consequence of  Theorem \ref{lemma:qs} and  Lemma \ref{lem:interlacing}.\newline

\begin{lemma}\label{lemma:reg:des0}
Consider a switching descriptor system described by
(\ref{eq:descriptor0})
%
and suppose that, for some $N_1\le N$,
 there is a symmetric positive-definite matrix $P$ satisfying
\begin{align}
PA_i^{-1}E_i + (A_i^{-1}E_i)^TP &< 0, \quad  i=1, \dots, N_1
 \label{eq:Q_1}\\
PA_j^{-1}E_j + (A_j^{-1}E_j)^TP &\le 0, \quad j=N_1+1, .., N
\end{align}
 and for each $j \in \{N_1+1, \cdots, N\}$ there is a subscript $i_j \in \{1, \cdots, N_1\}$ such that
{\small
\begin{equation}
\label{eq:rankCond3}
 rank(A_{i_j}^{-1}E_{i_j} -A_j^{-1}E_j) = rank(A_{i_j}^{-1}E_{i_j})- rank (A_j^{-1}E_j)
 \,.
\end{equation}
}
Also,
\begin{equation}
\label{eq:vdecrease}
x(t_*^+)Px(t_*^+) \le x(t_*^-)^TPx(t_*^-)
\end{equation}
if $\sigma$ switches at $t_*$.
 Then, the system is {GUES}.
 \end{lemma}
 
\vspace{.5em}
\noindent{\sc Proof:} We prove this result by showing that the hypotheses
 of Theorem \ref{lemma:qs} hold.
Since $P$ is positive definite,  hypothesis (a) holds.
Also, \eqref{eq:vdecrease} implies that   hypothesis (c) holds.
To see that hypothesis (b) holds, consider any $j \in \{N_1+1, \cdots, N\}$ and
  apply  Lemma \ref{lem:interlacing} with
$A=-PA_{i_j}^{-1}E_{i_j}$ and $B=-PA_j^{-1}E_j$ to obtain that
the kernel of $Q_j:=-PA_j^{-1}E_j -(A_j^{-1}E_j)^TP$ is the same as that of $-PA_j^{-1}E_j$ which also equals the kernel of
$A_j^{-1}E_j$; thus $Q_j$  and $A_j^{-1}E_j$ have the same kernel.
 Since $Q_j \ge 0$ and the kernel of $A_j^{-1}E_j$ and $\mathcal{C}(E_j, A_j)$ intersect only at zero,
  we conclude that $Q_j$ is positive definite on  the consistency
space of $(E_j, A_j)$. Hypothesis (b) now follows by taking into account (\ref{eq:Q_1}).
It now follows from Theorem \ref{lemma:qs} that
the switched system (\ref{eq:descriptor0}) is GUES.
\textbf{Q.E.D.}\newline

\noindent Now we consider a special class of switched descriptor systems described by
\begin{equation}\label{mixed}
E_{\sigma(t)} \dot x = A_{\sigma(t)}x\,, \qquad \sigma(t) \in \{  1, 2\}
\end{equation}
where each constituent system is stable,  with the first being index zero (standard system) and the second  index one;
also the rank of $A_1^{-1}E_1 -A_2^{-1}E_2$ is one.
We show that if the matrix $A_1^{-1}E_1A_2^{-1}E_2$ has no negative real eigenvalues, exactly one eigenvalue
at zero and some other regularity conditions hold then, the system
is GUES. To achieve this result, we recall the following result from \cite{shorten2009quadratic}.\newline

	\begin{theorem}\cite{shorten2009quadratic}
\label{th:marginally}
Suppose that $A$ is Hurwitz and all eigenvalues of $A-gh^T$
have negative real part, except one, which is zero.
Suppose also that
 $(A, g)$ is controllable and $(A, h)$ is observable.
 Then, there
exists a  matrix $P=P^T>0$ such that
\begin{align}
	A^TP+PA &< 0 \label{eq:lyap0a}\\
	(A-gh^T)^TP+P(A-gh^T) &\leq 0 \label{eq:lyap0b}
	\end{align}
if and only if the matrix product $A(A-gh^T)$ has no
real negative eigenvalues  and exactly one zero eigenvalue.
\end{theorem}
\vspace{.5em}

\noindent The following result follows from Lemma \ref{lemma:reg:des0} and Theorem \ref{th:marginally}.
%
%
%
%
%
%
%

\vspace{.5em}

\begin{theorem}\label{th:reg:des}
Consider a switching descriptor system described by
(\ref{mixed})
where
$x(\cdot)$ is continuous during switching and
suppose  that it satisfies the following conditions
  \begin{itemize}
  \item[(a)] $(E_1,A_1)$ and $(E_2,A_2)$ are   stable.
  \item[(b)] $(E_1,A_1)$ is  index zero and
   $(E_2,A_2)$ is index one.
   \item[(c)]\label{eq:speconetwo}
  There exists column matrices $g$ and $h$ such that
  \begin{equation}
  A_2^{-1}E_2 =A_1^{-1}E_1-gh^T,
  \end{equation}
 where $(A_1^{-1}E_1,\ g)$, $(A_1^{-1}E_1,\ h)$ are controllable and observable, respectively.
   \item[(d)] The matrix  $A_1^{-1}E_1A_2^{-1}E_2$ has no negative real eigenvalues and  exactly one zero eigenvalue.

\end{itemize}

\noindent Then the switching descriptor system (\ref{mixed}) is globally uniformly
exponentially stable about zero.
\end{theorem}
\vspace{5mm}
\noindent{\sc Proof:} We first show that the hypotheses of Theorem \ref{th:marginally} hold with $A = A_1^{-1}E_1$.
For $i=1,2$,  $(E_i, A_i)$ is stable; hence the non-zero eigenvalues of $A_i^{-1}E_i$ have negative real parts.
Since $(A_1, E_1)$ is index zero, $A_1^{-1}E_1$ is nonsingular and has no zero eigenvalues.
This implies that $A_1^{-1}E_1$ is Hurwitz.

Since $A_1^{-1}E_1A_2^{-1}E_2$ has  exactly one eigenvalue at zero, its nullity is one; the non-singularity of
$A_1^{-1}E_1A_2^{-1}$ now implies that the nullity of $E_2$ is one; hence the rank of $E_2$ and $A_2^{-1}E_2$ is  $n-1$.
Since  $(E_2, A_2)$ has index one and the nullity of $E_2$ is one, the matrix $A_2^{-1}E_2$ has a single eigenvalue at zero.
Thus,
all  eigenvalues of $A_2^{-1}E_2$ have negative real part except one which is zero.

 Recalling hypotheses (c) and (d) of  this theorem,
we see that the hypotheses of
Theorem \ref{th:marginally} hold with $A = A_1^{-1}E_1$. Hence there exists a matrix $P=P^T >0$ such that
\begin{eqnarray}
PA_1^{-1}E_1 + (A_1^{-1}E_1)^TP &<&0, \label{eq:Q_2}\\
 PA_2^{-1}E_2 + (A_2^{-1}E_2)^TP &\le&0.
\end{eqnarray}

\noindent Since $rank(A_1^{-1}E_1 - A_2^{-1}E_2) = rank(gh^T) = 1 =rank(A_1^{-1}E_1) -rank(A_2^{-1}E_2)$,
Lemma \ref{lemma:reg:des0} now implies that the switched descriptor system (\ref{mixed})
is globally uniformly exponentially stable about zero. \textbf{Q.E.D.}\newline

\begin{comment}
{\rm
The above result requires $x(\cdot)$ to be continuous during switching.
Since  $(E_1, A_1)$ is an index zero system   
  its consistency space is the whole state space.
Hence switching to this system can occur at any state.  
Since  $(E_2, A_2)$ is  index one,  the consistency space of this system is not the whole state space.
Thus, the switched system does not switch to the second system from
an arbitrary point in the state space. To switch to the second system, the state must be in the consistency space
of that system, that is it must be in $Im(A_2^{-1}E_2)$.
}
\end{comment}

\subsection* {Switching between index-one and index-two systems}
\noindent 
Now we consider switching between index-one  and   index-two descriptor systems. Our results in this subsection are based on an order reduction result from \cite{dimensionality}. 
They result from  an application of {\it full rank decomposition} to  a switched descriptor system in the form  of \eqref{eq:descriptor0}.\newline

\noindent{\it Full rank decomposition:} A pair of matrices
$(X,Y)$ is a {\it  decomposition} of $E\in \mathbb{R}^{n\times n}$ if
\begin{equation}
E=XY^T\,.
\end{equation}
If, in addition, $X$ and $Y$
both have full column rank we say that
 $(X,Y)$ is a {\it full rank decomposition} of $E$.
Note that, if $(X,Y)$ is a  full rank decomposition of $E\in \mathbb{R}^{n\times n}$ and
$rank(E) =r$ then, $X, Y \in\mathbb{R}^{n \times r}$ and   $rank(X) = rank(Y) = r$.  
Suppose $E$ has rank $r>0$ and  $\tilde{E} = Y^TA^{-1}X$   where $(X, Y)$ is a full rank decompostion of $E$;
then $(\tilde{E} , I)$ is a reduced order descriptor system with $r<n$ state variables. The original descriptor system $(E, A)$ is stable if and only if the non-zero eigenvalues of $\tilde{E}$ have negative real parts \cite{dimensionality}. 
Also, if  $k^*$ is the index of $(E, A)$  then  the index of the equivalent reduced order system $(\tilde{E},\, I)$ is $k^*\!-\!1$ \cite{dimensionality}. \newline

\noindent One can iteratively apply full rank decomposition to achieve further order reduction of $(\tilde{E} , I)$, provided that there
is  a decomposition $(\tilde{X}, \tilde{Y})$ of $\tilde{E}$ with $ \tilde{Y} \in \mathbb{R}^{r \times \tilde{r}}$
and $\tilde{r} <r$.
Since a non-zero  square matrix always has a full rank decomposition, one can always iteratively reduce a single linear system
$(E, A)$ to a standard system or to a system of algebraic equations, that is a system whose "E-matrix" is zero. \newline

Commonly, the switching condition on the state can be described by:
\begin{equation}
\label{eq:switchcond1}
x(t_*^+) = M_{ji}x(t_*^-)
\end{equation}
when $\sigma$ switches from $i$ to $j$ at $t_*$. 
Also, switching may be restricted in the sense that one does not switch from $i$ to $j$ at any state $x(t_*^-)$ in
$\mathcal{C}(E_i, A_i)$. In this case, the restriction may be described by
\begin{equation}
\label{eq:switchcond2}
x(t_*^-)\in \mathcal{S}_{ji}
\end{equation}

\begin{theorem}[Order reduction of linear switching descriptor systems  \cite{dimensionality}]
\label{lemma:red_order}
Consider a switching descriptor system described by
(\ref{eq:descriptor0})
 and switching conditions (\ref{eq:switchcond1}),(\ref{eq:switchcond2})
 when $\sigma$ switches from $i$ to $j$
and suppose that $(X_i, Y_i)$ is a  decomposition  of $E_i$  with $Y_i \in \mathbb{R}^{n \times r}$ for $i=1, \cdots N$.
%
Then, there exist matrices $T_1, \cdots, T_N$ such that  the following holds.
A function $x(\cdot)$ is a solution to  system (\ref{eq:descriptor0}) with (\ref{eq:switchcond1}),(\ref{eq:switchcond2}) if and only if
\begin{equation}
x(t)= T_{\sigma(t)} z(t)
\end{equation}
for all $t$
where $z(\cdot)$ is a solution to the
descriptor system
\begin{equation}
\label{eq:z system}
\tilde{E}_{\sigma(t)} \dot{z} =  z
\end{equation}
with
switching conditions
\begin{eqnarray}
z(t_*^+) &=& Y_j^TM_{ji} T_i z(t_*^-)         \label{eq:switchcond1z}\\
 T_iz(t_*^-) &\in& \mathcal{S}_{ji}        \label{eq:switchcond2z}
\end{eqnarray}
when $\sigma$ switches from $i$ to $j$ where
\begin{equation}
\tilde{E}_i = Y_i^TA_i^{-1}X_i \,.
\end{equation}
Moreover
\begin{equation}
z(t) = Y_{\sigma(t)}^T x(t)
\end{equation}
for all $t$,
$\mathcal{C}(\tilde{E}_i, I) = Y_i^T\mathcal{C}(E_i, A_i)$
and $z(\cdot)$ is continuous during switching if and only if the same is true of
$Y_{\sigma}^Tx$.
\ignore{
Then the behavior of $z$
is governed by the new switching descriptor system
\begin{eqnarray}
\label{eq:z system}
\tilde{E}_{\sigma(t)} \dot{z} = z\\
\label{eq:x and z}
x=T_i z
\end{eqnarray}
where
$
\tilde{E}_i = Y^TA_i^{-1}X_i$.
}
%
Hence,  
global uniform exponential stability (GUES) 
of the new system (\ref{eq:z system})-(\ref{eq:switchcond2z})
and the original system (\ref{eq:descriptor0})-(\ref{eq:switchcond2}) are equivalent.\newline
\end{theorem}

Now we present a general result which is a corollary of Theorem \ref{lemma:red_order} and Lemma \ref{lemma:reg:des0}.
\vspace{.5em}

\begin{corollary}\label{reg:des12}
Consider a switching descriptor system described by
(\ref{eq:descriptor0})
where
$Y_{\sigma}^Tx$ is continuous during switching
and
$(X_i, Y_i)$ is a decomposition  of $E_i$ with $Y_i \in \mathbb{R}^{n \times r}$ for $i=1, \dots, N$.
 Suppose that, for some $N_1\le N$, there is  a symmetric positive-definite matrix $P$ such that the
following conditions are satisfied, where $\tilde{E}_i =Y_i^TA_i^{-1}X_i$.
\begin{align}
P\tilde{E}_i + \tilde{E}_i^TP &<  0\,,\qquad  i=1, \dots, N_1
\label{eq:Q_i2} \\
P\tilde{E}_j + \tilde{E}_j^TP &\le 0\,, \qquad j=N_1+1, \dots, N
\end{align}
 and for each $j \in \{N_1+1, \dots, N\}$ there is a subscript $i_j \in \{1, \dots, N_1\}$ such that
\begin{equation}
 rank(\tilde{E}_{i_j} -\tilde{E}_j) = rank(\tilde{E}_{i_j})- rank (\tilde{E}_j)
 \,.
\end{equation}

\noindent Then, the system is globally uniformly exponentially stable about zero.
\end{corollary}
\vspace{.5em}

\noindent The above result  requires that $Y_{\sigma}^Tx$  be continuous during switching. 
Continuity of $Y_{\sigma}^Tx$   during switching  is equivalent to the following
switching condition. If  $\sigma$ switches from $i$ to $j$ at a point of discontinuity $t_*$ then,
\[
  Y_j^Tx(t_*^{+}) = Y_i^Tx(t_*^-) .
\]
Since $x(t_*^+)$ must be in $\mathcal{C}_j =\mathcal{C}(E_j, A_j)$, the above switch can only occur at states $x(t_*^-)$ in
$\mathcal{C}_i =\mathcal{C}(E_i, A_i)$ for which
\begin{equation}
\label{eq:sc2}
Y_i^Tx(t_*^-) \in Y_j^T\mathcal{C}_j\,.
\end{equation}

If  $(E_j, A_j)$ is index-one and
$Y_j \in \mathbb{R}^{r \times n}$ is full column rank where $r=rank(E)$ then,
switching to this system can occur from any state.
To see this, recall that the  kernel of $Y_j^T$  and $\mathcal{C}_j$ intersect only at the origin, and
since the system $(E_j, A_j)$ is index one,  the dimension of $\mathcal{C}_j$ is
$r=rank(E_j)$.  Hence the dimension of $Y^T\mathcal{C}_j$ is $r$.
Since  $Y_j^T \in \mathbb{R}^{r\times n}$ we now see that
$Y_j^T\mathcal{C}_j = \mathbb{R}^r$; hence (\ref{eq:sc2}) is satisfied for any $x(t_*^-)$.
This means that switching to  an index one system can occur from any state. 
For  an index-two system
\begin{eqnarray}
 Y^T_i\mathcal{C}_i&=&\text{Im}(Y^T_i(A_i^{-1}E_i)^2)\nonumber\\
 &=&\text{Im}((Y^T_iA_i^{-1}X_i)^2Y^T_i)\nonumber\\
 &=&\text{Im}((Y^T_iA_i^{-1}X_i)^2)\,.
\end{eqnarray}
\noindent 
For an index two system, $Y^T_iA_i^{-1}X_i$ is singular; hence the dimension of  $ Y^T_i\mathcal{C}_i$ is stricltly less than $r$.
 Hence we can always find $x(t_*^-) $ such that $Y_i^Tx(t_*^-) \notin Y_j^T\mathcal{C}_j$. Thus we cannot arbitrarily switch to an index-two system. 

\noindent Now to conclude we present our next main result:  switching between an index-one  and   an index-two descriptor system. The following result is obtained from Corollary \ref{reg:des12} and Theorem \ref{th:marginally}.\newline

\begin{theorem}\label{index12} Consider a switching descriptor system described by
\begin{equation}
\label{dswisys12}
E_{\sigma(t)} \dot x = A_{\sigma(t)}x\,, \qquad \sigma(t) \in \{1, 2\},
\end{equation}
where $Y_{\sigma}^Tx$ is continuous during switching
and
$(X_i, Y_i)$ is a full rank decomposition  of $E_i$ with $ Y_i \in \mathbb{R}^{n \times r}$ for $i=1,2$.
 Suppose that the
following conditions are satisfied where $\tilde{E}_i =Y_i^TA_i^{-1}X_i$ for $i=1, 2$.
  \begin{itemize}
  \item[(a)]$(E_1,A_1)$ and $(E_2,A_2)$ are   stable.
  \item[(b)] $(E_1,A_1)$ is  index one and
   $(E_2,A_2)$ is index two.
   \item[(c)]
  There exist vectors  $g$ and $h$ such that
  \begin{equation}
  \label{eq:bc2}
  \tilde{E}_2 = \tilde{E}_1-gh^T
  \end{equation}
with $(\tilde{E}_1, g)$ controllable and $(\tilde{E}_1, h)$ observable.
   \item[(d)]
    $\tilde{E}_1\tilde{E}_2$
has no negative real eigenvalues and  exactly one zero eigenvalue.

\end{itemize}
\noindent Then the switched descriptor system (\ref{dswisys12}) is globally uniformly
exponentially stable about zero.
\end{theorem}

\vspace{5mm}

\noindent{\sc Proof:}
Since $(X_1, Y_1)$ is a full rank decomposition of $E_1$ and $(A_1, E_1)$ is stable  and
index-one, its corresponding reduced order system $(\tilde{E}_1, I)$ is stable and index-zero.
Since $(X_2, Y_2)$ is a full rank decomposition of $E_2$ and $(A_2, E_2)$ is stable  and
index-two, its corresponding reduced order system $(\tilde{E}_2, I)$ is stable and index one.
Theorem \ref{th:reg:des} now guarantees {GUES} of
the reduced-order switched system.
Theorem \ref{lemma:red_order} now implies the same stability properties for
the original switched system (\ref{dswisys12}).
\textbf{Q.E.D.}\newline

\begin{comment}
\noindent Clearly the continuity assumption on $Y_{\sigma}^Tx$ restricts the applicability of our results to
certain decompositions of $E_i$, since full rank decompositions are in general non-unique.
Note however that, for any system,  if $E_i = E$ for all $i$  and $Ex$ is continuous then $Y^Tx$ is continuous
for any  full rank decomposition   $(X, Y)$  of $E$.
\end{comment}
%

\section{Numerical examples}

\begin{figure}[h]
\centering
\includegraphics[scale=0.32]{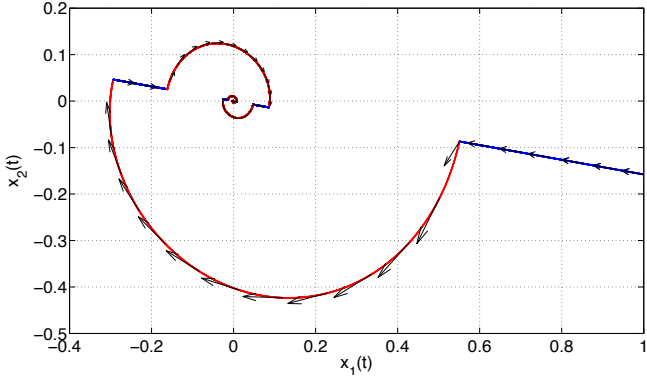}
\caption{Switching between an index-zero and  an index-one descriptor system.}
\label{fig:switchingZeroOne}
\end{figure}

\begin{example}[Switching between an index-zero and an index-one descriptor system]
\label{ex:switch0_1}
Consider a switched system of the form
\eqref{mixed}
where  $x(\cdot)$ is continuous and
\begin{align*}
&E_{1}= \begin{bmatrix}
1 & 0  \\
0 & 1 
\end{bmatrix},\,
\qquad
A_{1}= \begin{bmatrix}
-1 & 4\pi  \\
-4\pi & -4 
\end{bmatrix}&
\\
&E_2=\begin{bmatrix}
1 &  0\\
1 & 0
\end{bmatrix},
\qquad
A_2= \begin{bmatrix}
-1 & 4\pi \\
-\pi & -1
\end{bmatrix}
&
\end{align*}
\noindent Note that $(E_{1},A_{1})$ is a stable  index-zero system
whereas $(E_2, A_2)$ is a stable index one descriptor system
 whose consistency space $\mathcal{C}_2$ is  $\text{Im}\left(\begin{bmatrix} 1& -k_1k_2\end{bmatrix}^T\right)$ 
 where $k_1=\pi-1$ and $k_2=1/(4\pi +1)/$;
 note that  $\mathcal{C}_2$ can be represented by the line $(k_1k_2) x_1+x_2=0$. Note also that
$A^{-1}_1E_{1}-A^{-1}_2E_2=gh^T$,
where $g^T=\begin{bmatrix}1 & {1}/{4\pi}\end{bmatrix}$ and $h=\begin{bmatrix}-{4\pi}/{k_3} & {\pi}/{k_3}\end{bmatrix}^T$ with $k_3=4\pi^2+1$.
The pairs $(A^{-1}_1E_{1},g)$ and $(A^{-1}_1E_{1},h)$ are controllable and observable, respectively. The eigenvalues
 of $A^{-1}_1E_{1}A^{-1}_2E_2$ are  $( 0, 0.0042)$.
  Hence from Theorem \ref{th:reg:des}, the switched system described above is globally uniformly exponentially stable about zero.\newline 
  
 \noindent To illustrate GUES of  this system we consider a special switching signal. 
 The switching signal cannot be arbitrary, because of the assumption that   $x(\cdot)$ is continuous during switching. 
 When switching from the index-zero system $(E_{1},A_{1})$ to the index-one system $(E_2,A_2)$ at a time $t_*$, 
  we  must have  $x(t_*^-) \in \mathcal{C}_2$.  
 However, switching from the  index-one system to the index-zero can happen at any arbitrary time.
 
 \noindent The restriction to switch only when consistency spaces intersect  can be enforced through state dependent switching. However, for the purpose of illustration we consider a switching signal which is combination of state dependent switching and periodic switching. 
 To explain, let $t_*$ be a time  when the trajectory of the index-zero system  reaches $\mathcal{C}_2$, i.e. $x(t_*) \in \mathcal{C}_2$ or equivalently $(k_1k_2) x_1(t_*)+x_2(t_*)=0$. Now, we let   \eqref{mixed} switch from  $(E_{1},A_{1})$ to $(E_2,A_2)$, i.e.,  
$\sigma(t_*^-)= 1$ and $\sigma(t_*^+)= 2$. Every time this switch happens we fix $\sigma(t)=2$ for a time period $T$ before the system switches back to $(E_{1},A_{1})$.\newline


\noindent  Now we plot the trajectory of \eqref{mixed} with $T=0.2$ seconds and the initial state 
$x_0=\begin{bmatrix}1 & -k_1k_2\end{bmatrix}^T=\begin{bmatrix}1 & -0.1579\end{bmatrix}^T \in \mathcal{C}(E_2,A_2)$ using MATLAB (code is available online at \cite{corless_code}). The resulting trajectory is illustrated in Figure \ref{fig:switchingZeroOne} for $0 \leq t\leq 2$ seconds. 
The trajectory for the index-zero system is represented by the red spiral and the trajectory for the index-one system is along the blue line passing through origin.
\end{example}

\begin{example}[Switching between an index-one and an index-two descriptor system]

\noindent Consider a switched system of the form
\eqref{dswisys12} where

\begin{align*}
&E_1= \begin{bmatrix}
1 & 0  & 0\\
0 & 0 & 0\\
0 & 0 & 1
\end{bmatrix},\qquad
A_1= \begin{bmatrix}
-1 & 0 & 4\pi  \\
0 & -1 & 0\\
-4\pi &0 & -4 
\end{bmatrix}&
\\
&E_2=\begin{bmatrix}
1 & 0  & 0\\
0 & 0 & 1\\
0 & 0 & 1
\end{bmatrix} ,
\qquad
A_2= \begin{bmatrix}
-k_2k_3 & 0 &0 \\
0 & -1 &0\\
-4k_1k_2k_3 & -1 &-4k_3
\end{bmatrix}
&
\end{align*}
\noindent 
with $k_1, k_2$ and $k_3$ as defined in the previous example.
Note that $(E_1,A_1)$ and $(E_2, A_2)$ are a stable index one and index two descriptor systems  respectively. 
A full rank decomposition $(X_1, Y_1)$ of $E_1$ is given by
 \[
 X_1= 
 Y_1=\begin{bmatrix}1 &0 &0\\0 & 0 &1\end{bmatrix}^T
 \]
and  a full rank decomposition of $(X_2, Y_2)$ of $E_2$ is given by
\[
X_2=\begin{bmatrix}1 &0 & 0 \\0 & 1 &1\end{bmatrix}^T  , \qquad 
Y_2= Y_1 \,.
\]  Now we use Theorem  \ref{lemma:red_order} to obtain the equivalent reduced order switched system
\begin{eqnarray}\label{numerical_eg2}
\tilde{E}_{\sigma(t)} \dot{z}(t) = z(t),
\end{eqnarray}
where  $\tilde{E}_1=Y^T_1{A}^{-1}_1X_1$,  $\tilde{E}_2=Y^T_2{A}^{-1}_2X_2$ and $z(t)=Y^Tx(t)$  with $Y=Y_1=Y_2$.  
Upon evaluating $\tilde{E}_1$ and $\tilde{E}_2$ we can observe that \eqref{numerical_eg2} is the same as the switched descriptor system described in Example \ref{ex:switch0_1}. Now, if we assume that $z(\cdot)$ is continuous during switching then  it follows from the conclusions in   Example \ref{ex:switch0_1} and Theorem \ref{lemma:red_order} that \eqref{numerical_eg2} is GUES. 
One can also use  Theorem \ref{index12} to deduce GUES.
\end{example}

\section{Conclusions}

In this paper we derive stability conditions for a switched system where  switching occurs between linear  descriptor systems  of non-homogeneous indices.
To the best of our knowledge, these conditions are some of the first to consider the case of switching between modes of different indices. 
For specific cases, such as switching between  two systems whose indices differ by  one,  spectral conditions are derived that can be 
used to check   stability of such systems in an elementary manner.  Examples are also given to illustrate the use of our results. Future work will 
consider extending our analysis using non-quadratic Lyapunov functions and also consider  switching between  systems whose indices differ by more than one.

\end{document}